\documentclass[twoside,11pt]{article}
\usepackage[english]{babel}
\usepackage{amsmath,amssymb,amsfonts,bm,amsthm}
\usepackage{cite}
\usepackage{indentfirst}
\usepackage{ifpdf}
\usepackage{graphicx}

\ifpdf
\usepackage[bookmarks=false,pdfstartview=FitH,linkbordercolor={0.5 1 1},%
citebordercolor={0.5 1 0.5},unicode,%
pagebackref,hyperindex]{hyperref}
\fi
\usepackage{url}

\newtheoremstyle{myplain}
{4pt}
{4pt}
{\itshape}
{\parindent}
{\bfseries}
{.}
{.5em}
{}

\newtheoremstyle{myremark}
{4pt}
{4pt}
{\rmfamily}
{\parindent}
{\bfseries}
{.}
{.5em}
{}

\theoremstyle{myplain}
\newtheorem{theorem}{Theorem}
\newtheorem{lemma}{Lemma}
\newtheorem*{corollary*}{Corollary}

\theoremstyle{myremark}

\newtheorem{example}{Example}

\newlength{\ArticleHeadDepth}
\newlength{\captwidth}%
\makeatletter
\def\@listI{\leftmargin\leftmargini
            \parsep 0\p@ \@plus1\p@ \@minus\p@
            \topsep 8\p@ \@plus2\p@ \@minus4\p@
            \itemsep0\p@}
\let\@listi\@listI
\@listi
\def\@listii {\leftmargin\leftmarginii
              \labelwidth\leftmarginii
              \advance\labelwidth-\labelsep
              \topsep    0\p@ \@plus2\p@ \@minus\p@}
\def\@listiii{\leftmargin\leftmarginiii
              \labelwidth\leftmarginiii
              \advance\labelwidth-\labelsep
              \topsep    0\p@ \@plus\p@\@minus\p@
              \parsep    \z@
              \partopsep \p@ \@plus\z@ \@minus\p@}

\renewcommand\labelitemii{$\m@th\bullet$}

\def\sup{\mathop{\operator@font sup}}
\def\inf{\mathop{\operator@font \vphantom{sup}inf}}
\makeatother

\makeatletter
\renewenvironment{thebibliography}[1]
     {\section*{\refname}
      \def\@biblabel##1{##1.}
      \small
      \list{\@biblabel{\@arabic\c@enumiv}}%
           {\settowidth\labelwidth{\@biblabel{#1}}%
            \leftmargin\labelwidth
            \advance\leftmargin\labelsep
            \usecounter{enumiv}%
            \let\p@enumiv\@empty
            \renewcommand\theenumiv{\@arabic\c@enumiv}}%
            \renewcommand\newblock{\hskip .11em \@plus.33em \@minus.07em}%
      \sloppy\clubpenalty4000\widowpenalty4000%
      \sfcode`\.=\@m}
     {\def\@noitemerr
       {\@latex@warning{Empty `thebibliography' environment}}%
      \endlist}

\def\@maketitle{\markboth{}{}%
  \center
   {\LARGE\bfseries\boldmath
  \pretolerance=10000
  \settoheight{\ArticleHeadDepth}{X}
  \addtolength{\ArticleHeadDepth}{-\headsep}
  \addtolength{\ArticleHeadDepth}{32.4pt}
  \mbox{}\vskip\ArticleHeadDepth
  \@title
  \par}\vskip 10pt
 {\lineskip .5em\noindent\ignorespaces {\large\bfseries\@author}\vskip8pt}
 {\ignorespaces
 \centering\itshape\small\@institute\par}
 }

\renewenvironment{abstract}
      {\list{}{\setlength\leftmargin{2\leftmargini}\rightmargin\leftmargin\small
      \labelwidth=\z@
      \listparindent=\z@
      \labelsep=\z@
      \itemindent\listparindent}%
                \item[\hskip\labelsep
                     \bfseries\abstractname]}
               {\endlist}

\renewcommand\section{\@startsection{section}{1}{\z@}%
                       {-18\p@ \@plus -4\p@ \@minus -4\p@}%
                       {4\p@ \@plus 2\p@ \@minus 1\p@}%
                       {\normalfont\normalsize\center\MakeUppercase
                        }}

\def\appendix#1{\par\setcounter{section}{0}\setcounter{subsection}{0}
\setcounter{equation}{0} \vskip 18\p@ \@plus 4\p@
\@minus 4\p@
\begin{flushright}{\slshape\MakeUppercase{\appendixname}~#1}\end{flushright}
\vskip 4\p@ \@plus 2\p@ \@minus 1\p@
\def\thesection{A.\arabic{section}}
\def\theequation{A.\arabic{equation}}
}

\def\institute#1{\gdef\@institute{#1}}
\def\received#1{\gdef\@received{#1}}\edef\@received{When received?}
\makeatother

\begin{document}

\title{Algebraic Unsolvability of Problem of Absolute Stability of
Desynchronized Systems Revisited}

\author{V. S.~Kozyakin}

 \institute{Institute for Information Transmission Problems RAS, Moscow, Russia}

\maketitle

\begin{abstract}~
In the author's article ``Algebraic unsolvability of problem of absolute
stability of desynchronized systems'' (Automat. Remote Control 51 (1990),
no. 6, pp. 754--759), it was shown that in general for linear
desynchronized systems there are no algebraic criteria of absolute
stability. In this paper, a few misprints occurred in the original version
of the article are corrected, and two figures are added.
\end{abstract}

\section{Introduction}

In complex control systems containing sampled-data elements, it is possible
that these elements operate asynchronously. In some cases asynchronous
character of operation of sampled-data elements does not influence stability
of system. In other cases any small desynchronization of the updating moments
of sampled-data elements leads to dramatic changes of dynamics of a control
system, and the system loses stability \cite{KKKK:AiT83-84:e}. Last years
there is begun (see, e.g.,
\cite{KKKK:AiT83-84:e,Klep83:e,KKKK:DAN84:e,KKKK:MCS84,Klep:AIT85:e})
intensive studying of the effects connected with asynchronous operation of
control systems; both necessary, and sufficient stability conditions for
various classes of asynchronous systems were obtained. At the same time no
one succeed in finding general, effectively verified criteria of stability of
asynchronous systems, similar to known for synchronous systems
\cite{Zsyp:63:e}. The problem on stability of linear asynchronous systems has
appeared more difficult than the problem on stability of synchronous systems.
In the paper, attempt of formal explanation of complexity of the stability
analysis problem for linear asynchronous systems is undertaken. It is shown
that there are no criteria of absolute stability of linear asynchronous
systems consisting of a finite number of arithmetic operations.

\section{Statement of Problem}\label{S-problem}

Consider a discrete-time linear control system whose dynamics is described by
the vector difference equation
\begin{equation}\label{E:main}
x(n)=A(n)x(n-1)\quad	(n= 1,2,\ldots),
\end{equation}
where $x(n)=\{x_{1}(n), x_{2}(n), \ldots , x_{N}(n)\}$ is the state vector of
the system and $A(n)=(a_{ij}(n))$ is a square matrix of dimension $N$ with
the elements $a_{ij}(n)$.

The system \eqref{E:main} will be called synchronous if
$A(n)\equiv\mathrm{const}$. If $A(n)\not\equiv\mathrm{const}$, and the set
$\{A(n): n=1,2,\ldots\}$ consists of finitely many elements $A_{1},
A_{2},\ldots,A_{M}$, then the system \eqref{E:main} will be called
asynchronous or desynchronized.

Let $\mathfrak{A}=\{A_{1}, A_{2},\ldots,A_{M}\}$ be a finite totality of
square matrices of dimension $N$. The system \eqref{E:main} will be called
absolutely stable with respect to the class of matrices $\mathfrak{A}$ (cf.
\cite{AizGant:e}) if there exists $c=c(\mathfrak{A})$ such that for any
sequence of matrices $A(n)\in\mathfrak{A}$ the following estimates hold:
\begin{equation}\label{E:estim}
\|A(n)A(n-1)\cdots A(1)x\| \le c\|x\|\quad (n=1,2,\ldots).
\end{equation}

Let us call the system \eqref{E:main} absolutely exponentially stable with
respect to the class of matrices $\mathfrak{A}$ if there exist
$c=c(\mathfrak{A})$ and $q=q(\mathfrak{A})<1$ such that for any sequence of
matrices $A(n)\in\mathfrak{A}$ the following estimates hold:
\begin{equation}\label{E:estim-exp}
\|A(n)A(n-1)\cdots A(1)x\| \le cq^{n}\|x\|\quad (n=1,2,\ldots).
\end{equation}

If the class of matrices $\mathfrak{A}$ consists of the square matrices
$A_{1}=(a_{1ij})$, $A_{2}=(a_{2ij})$, \ldots, $A_{M}=(a_{Mij})$ of dimension
$N$ then, for its description, it suffices to specify $MN^{2}$ numbers:
$a_{111}$, $a_{112}$, \ldots, $a_{1NN}$, $a_{211}$, $a_{212}$, \ldots ,
$a_{2NN}$, \ldots, $a_{M11}$, $a_{M12}$, \ldots, $a_{MNN}$. Therefore, each
class $\mathfrak{A}$ consisting of $M$ square matrices of dimension $N$ can
be treated as a point in some space $\mathfrak{M}(M,N)=R^{MN^{2}}$. Denote by
$S(M,N)$ the set of those classes $\mathfrak{A}$ in the space
$\mathfrak{M}(M,N)$ with respect to which the system \eqref{E:main} is
absolutely stable. By $E(M,N)$ we denote the set of those classes
$\mathfrak{A}\in\mathfrak{M}(M,N)$ with respect to which the system
\eqref{E:main} is absolutely exponentially stable.

Now, the problem of studying the absolute stability of the system
\eqref{E:main} can be reformulated as the problem of description of the sets
$S(M,N)$ and $E(M,N)$; the simpler in some sense the structure of the sets
$S(M,N)$ or $E(M,N)$ the easier to obtain a criterion of absolute stability
or absolute exponential stability.

The sets $S(1,N)$ and $E(1,N)$ allow a simple description. Indeed, each class
$\mathfrak{A}\in\mathfrak{M}(1,N)$ consists of a single matrix. Therefore, we
need to obtain conditions of stability or asymptotic stability of some
difference equation $x(n)=Ax(n-1)$. The Routh--Hurwitz stability criterion
\cite{Zsyp:63:e} allows to represent these conditions as a finite system of
polynomial inequalities including the elements $a_{ij}$ of the matrix $A$.
Verification of the obtained inequalities can be performed by a finite number
of arithmetic operations over the elements of the matrix $A$. In other words,
the question whether an arbitrary class $\mathfrak{A}=\{A\}$ belongs to the
set $S(1,N)$ or $E(1,N)$ may be resolved by a finite number of arithmetic
operations.

Is it possible, for $M\ge 2$, by a finite number of arithmetic operations to
resolve the question whether an arbitrary class
$\mathfrak{A}\in\mathfrak{M}(M,N)$ belongs to the set $S(N,N)$ or $E(M,N)$?
The answer to this question will be given in the next section.

\section{Main Result}\label{S:main}

Let $u=\{u_{1}, u_{2},\ldots,u_{L}\}$ be an element of the coordinate space
$R^{L}$. A finite sum $p(u) = \sum p_{i_{1}i_{2}\ldots i_{L}} u_{1}^{i_{1}}
u_{2}^{i_{2}}\cdots u_{L}^{i_{L}}$ with numerical coefficients
$p_{i_{1}i_{2}\ldots i_{L}}$ is called a polynomial in variable $u\in R^{L}$.
A set $U\subseteq R^{L}$ is said to have the \textit{SA}-property
\cite{Trev:e} if there exists a finite number of polynomials $p_{1}(u)$,
\ldots , $p_{k}(u)$, $p_{k+1}(u)$, \ldots, $p_{k+l}(u)$ such that $U$
coincides with the set of elements $u\in R^{L}$ satisfying the condition
\begin{equation}\label{E:sa-set}
p_{1}(u)>0,\ldots , p_{k}(u)>0,\quad
p_{k+1}(u)=\ldots= p_{k+l}(u)=0.
\end{equation}

A set $U$ is called semialgebraic \cite{Trev:e} if it is a unity of a finite
number of the sets possessing the \textit{SA}-property.
\begin{theorem}\label{T:1}
Let $M,N\ge2$. If a subset $U$ of the space $\mathfrak{M}(M,N)$ satisfies
conditions $E(M,N)\subseteq U\subseteq S(M,N)$ then it is not semialgebraic.
\end{theorem}

The proof of the theorem is given in the Appendix.

Semialgebraicity of a set is equivalent to the existence of a criterion
(consisting in verification of a finite number of the conditions of the form
\eqref{E:sa-set}) which allows by a finite number of arithmetic operations of
addition, subtraction, multiplication and comparison of numbers to establish
belonging of an element to a given set. As seen from Theorem~\ref{T:1},
neither the set $S(M,N)$ nor the set $E(M,N)$ are semialgebraic. So, the
meaning of Theorem~\ref{T:1} is that in general, by a finite number of
arithmetic operations, it is impossible to ascertain whether a desynchronized
system \eqref{E:main} is absolutely stable (absolutely exponentially stable)
or not.

The problem on the existence of algebraic criteria of stability is acute also
for classes of desynchronized systems different from those considered above.
For example, in the theory of continuous-time desynchronized systems there
arises the problem of stability of the so-called regular systems
\cite{KKKK:DAN84:e,KKKK:MCS84} (i.e., the systems with the infinite number of
updating moments for each component). The discrete-time system
\eqref{E:main}, and the related to it sequence of matrices $A(n)$, will be
called regular if each matrix $A_{i}$ from the class $\mathfrak{A}=\{A_{1},
A_{2},\ldots,A_{M}\}$ appears in the sequence $\{A(n)\}$ infinitely many
times. Denote by $r(n)$ the greatest integer $r$ having the property: the set
of matrices $\{A(1),A(2),\ldots,A(n)\}$ can be decomposed in $r$ subsets
$\{A(1),\ldots,A(n_{1})\}$, $\{A(n_{1}+1),\ldots,A(n_{2})\}$, \ldots,
$\{A(n_{r-1}+1),\ldots,A(n)\}$ such that each of them contains all the
matrices $A_{1}, A_{2},\ldots,A_{M}$. Clearly, the system \eqref{E:main} is
regular if and only if $r(n)\to\infty$.

The system \eqref{E:main} will be called absolutely exponentially regularly
stable with respect to the class of matrices $\mathfrak{A}=\{A_{1},
A_{2},\ldots,A_{M}\}$ if there exist $c=c(\mathfrak{A})$,
$q=q(\mathfrak{A})<1$ such that for any regular sequence of matrices
$A(n)\in\mathfrak{A}$ the following estimates hold:
\[
\|A(n)A(n-1)\cdots A(1)x\| \le cq^{r(n)}\|x\|.
\]
Denote by $R(M,N)$ the set all the classes $\mathfrak{A}\in\mathfrak{M}(M,N)$
with respect to which the system \eqref{E:main} is absolutely exponentially
regularly stable.
\begin{theorem}\label{T:2}
Let $M,N\ge2$. Then the set $R(M,N)$ is not semialgebraic.
\end{theorem}

To prove Theorem~\ref{T:2} it suffices to note that the set $R(M,N)$ contains
$E(M,N)$, and is contained in $S(M,N)$. Then by Theorem~\ref{T:1} it is not
semialgebraic. In other words, for $M,N\ge2$ there are no semialgebraic
criteria of absolutely exponentially regular stability of discrete-time
desynchronized systems.

\section{Addendum}

As shown above, the problem of absolute stability of the system
\eqref{E:main} can be reduced to the analysis of behaviour of infinite
products of the matrices $A(n)\in\mathfrak{A}$. Theorem~\ref{T:3} below
reduces the same problem to the descriptive-geometric question on existence
in the space $R^{N}$ such a norm in which each matrix $A_{1},
A_{2},\ldots,A_{M}$ is contractive.

\begin{theorem}\label{T:3}
The system \eqref{E:main} is absolutely stable in a class of matrices
$\mathfrak{A}=\{A_{1}, A_{2},$ \ldots, $A_{M}\}$ if and only if there is a
norm $\|\cdot\|$ in $R^{N}$ for which the following inequalities hold:
\begin{equation}\label{E:ineq1}
\|A_{1}\|, \|A_{2}\|,\ldots, \|A_{M}\|\le1.
\end{equation}

The system \eqref{E:main} is absolutely exponentially stable in a class of
matrices $\mathfrak{A}$ if and only if there is a norm $\|\cdot\|$ in $R^{N}$
and a number $q<1$ for which the following inequalities hold:
\begin{equation}\label{E:ineqq}
\|A_{1}\|, \|A_{2}\|,\ldots, \|A_{M}\|\le q.
\end{equation}
\end{theorem}

The proof of the theorem is given in the Appendix.

Several important properties of the sets $S(M,N)$ and $E(M,N)$ follow from
Theorem~\ref{T:3}. For example, the set $E(M,N)$ is open in
$\mathfrak{M}(M,N)$; the set $E(M,N)$ belongs to the interior of the set
$S(M,N)$.

Due to openness of the set $E(M,N)$, if the system \eqref{E:main} is
absolutely exponentially stable with respect to some class
$\mathfrak{A}^{0}=\{A_{1}^{0}, A_{2}^{0},\ldots,A_{M}^{0}\}$ then it is also
absolutely exponentially stable with respect to any class
$\mathfrak{A}=\{A_{1}, A_{2}$, \ldots, $A_{M}\}$ of matrices $A_{i}$
sufficiently close to the corresponding matrices $A_{i}^{0}$
$(i=1,2,\ldots,M)$.

Theorems~\ref{T:1} and \ref{T:3} imply that the problem of construction, for
a given set of square matrices, of a norm satisfying conditions
\eqref{E:ineq1} or \eqref{E:ineqq} is algebraically unresolvable.

Theorem~\ref{T:1} states that in general there are no effective criteria of
absolute stability of desynchronized systems \eqref{E:main}. Nevertheless,
such criteria may exist for some particular desynchronized systems. Let us
present examples.

\begin{example}\label{Ex:1}
Denote by $\mathfrak{R}(N)$ the subset of the space $\mathfrak{M}(N,N)$
consisting of the classes $\mathfrak{A}=\{A_{1}, A_{2},\ldots,A_{N}\}$ of
matrices $A_{i}$ of the form
\begin{equation}\label{E:Amix}
A_{i}=\left(
\begin{array}{cccccc}
1&0&\ldots&0&\ldots&0\\
0&1&\ldots&0&\ldots&0\\
\ldots&\ldots&\ldots&\ldots&\ldots&\ldots\\
a_{i1}&a_{i2}&\ldots&a_{ii}&\ldots&a_{iN}\\
\ldots&\ldots&\ldots&\ldots&\ldots&\ldots\\
0&0&\ldots&0&\ldots&1
\end{array}\right).
\end{equation}

The problem on absolute stability of the system \eqref{E:main} with respect
to the classes $\mathfrak{A}\in\mathfrak{R}(N)$ arises in
\cite{KKKK:AiT83-84:e,Klep83:e,KKKK:DAN84:e,KKKK:MCS84,Klep:AIT85:e} in the
process of study of continuous-time systems with a special types of
desynchronization of updating moments.

\begin{theorem}\label{T:4}
The system \eqref{E:main} is absolutely stable with respect to the class of
matrices $\mathfrak{A}=\{A_{1},A_{2}\}\in\mathfrak{R}(2)$ of the form
\eqref{E:Amix} if and only if one of the following system of relations holds:
\begin{description}
\item[\quad\rm a)] $a_{11}=1$, $a_{12}=a_{21}=0$, $a_{22}=1$;
\item[\quad\rm b)] $a_{11}=1$, $a_{12}=0$, $a_{22}=-1$, $a_{21}$ is
    arbitrary;
\item[\quad\rm c)]$a_{11}=-1$, $a_{12}$ is arbitrary, $a_{21}=0$,
    $a_{22}=1$;
\item[\quad\rm d)] $a_{11}=a_{22}=-1$, $0\le a_{12},a_{21}<4$;
\item[\quad\rm e)] $|a_{11}|<1$, $|a_{22}|<1$,
    $-(1-|a_{11}|)(1-|a_{22}|)\le a_{12}a_{21}\le(1-a_{11})(1-a_{22})$.
\end{description}
\end{theorem}

The proof of Theorem~\ref{T:4} is cumbersome and so is skipped. Let us point
out that the criterion of absolute stability of the system \eqref{E:main}
with respect to the classes of matrices from $\mathfrak{R}(2)$, given in
Theorem~\ref{T:4}, is semialgebraic.
\end{example}

\begin{example}\label{Ex:2}
Denote by $\mathfrak{R}_{+}(N)$ the subset of the space	$\mathfrak{M}(N,N)$
consisting of the classes $\mathfrak{A}=\{A_{1}, A_{2},\ldots,A_{N}\}$ of
matrices $A_{i}$ of the form \eqref{E:Amix} for which $a_{ij}>0$. The set
$S(N,N)\cap \mathfrak{R}_{+}(N)$ is semialgebraic. The criterion of absolute
stability of the system \eqref{E:main} with respect to the classes of
matrices $\mathfrak{A}$ from $\mathfrak{R}_{+}(N)$ consists in verification
that the maximal eigenvalue of the matrix $A=(a_{ij})$ does not exceed $1$.
This assertion is proved similarly to Theorem~2 from the first part of
\cite{KKKK:AiT83-84:e}.
\end{example}

\appendix{}

1. \emph{Proof of Theorem~\ref{T:1}} suffices to present for the case
$M=N=2$. The idea of proof is simple. We construct two families of matrices
depending on the real parameter $t\in[-1,1]$:
\begin{equation}\label{E:A1}
G(t)=(1-t^{4})\left(\begin{array}{cc}
1& -\frac{t}{\sqrt{1-t^{2}}}\\
0&0
\end{array}\right),\quad
H(t)=(1-t^{4})\left(\begin{array}{cc}
1-2t^{2}& -2t\sqrt{1-t^{2}}\\
2t\sqrt{1-t^{2}}&1-2t^{2}
\end{array}\right).
\end{equation}

\begin{figure}[htbp!]
\centering\includegraphics*{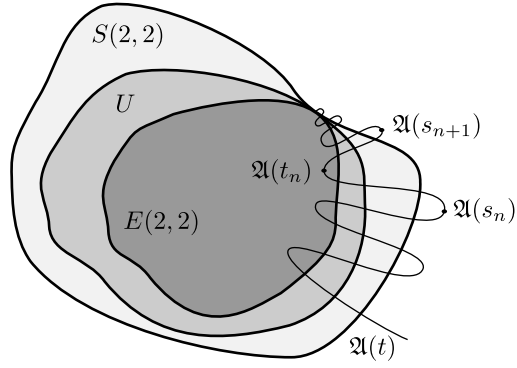}
\caption{A case when the set
$U\cap W$ has infinitely many components of connectedness}\label{semialg}
\end{figure}

The set $W$ of all the classes $\mathfrak{A}$ of the form
$\mathfrak{A}=\mathfrak{A}(t)= \{G(t), H(t)\}$ forms in the space
$\mathfrak{M}(2,2)$ an algebraic set. Suppose that the set $U$ is
semialgebraic. Then the set $U\cap W$ is also semialgebraic. Therefore, by
the theorem of Whitney (see, e.g., \cite{Milnor:e}) about a finite number of
connected components of a real algebraic set, a neighbourhood of the class
$\mathfrak{A}(0)$ in $U\cap W$ should be either empty or consisting of a
finite number of the connected components. We will show that it has
infinitely many components of connectedness, see Fig.~\ref{semialg}. So, the
set $U$ is not semialgebraic.

Let us pass to the theorem proof. Denote by $|\cdot|$ the Euclidean norm in
$R^{2}$: if $x=\{\xi,\eta\}$ then $|x|=\sqrt{\xi^{2}+\eta^{2}}$. Consider two
families of matrices:
\begin{equation}\label{E:A2}
P(\varphi)=\left(\begin{array}{cc}
1& -\tan\varphi\\
0&0
\end{array}\right),\quad
R(\varphi)=\left(\begin{array}{cr}
\cos 2\varphi& -\sin 2\varphi\\
\sin 2\varphi&\cos 2\varphi
\end{array}\right).
\end{equation}
\begin{lemma}\label{L:A1}
Let $\varphi=\pi/(2n+1)$. Then
$(PR^{n}P)(\varphi)=-(\cos\varphi)^{-1}P(\varphi)$.
\end{lemma}

\begin{lemma}\label{L:A2}
Let $\varphi=\pi/(2n)$. Then	$(PR^{m}P)(\varphi)=\lambda_{m,n}P(\varphi)$,
where $|\lambda_{m+n,n}|=|\lambda_{m,n}|$, $|\lambda_{m,n}|\le 1$.
\end{lemma}

\begin{figure}[htbp!]
\centering\includegraphics*{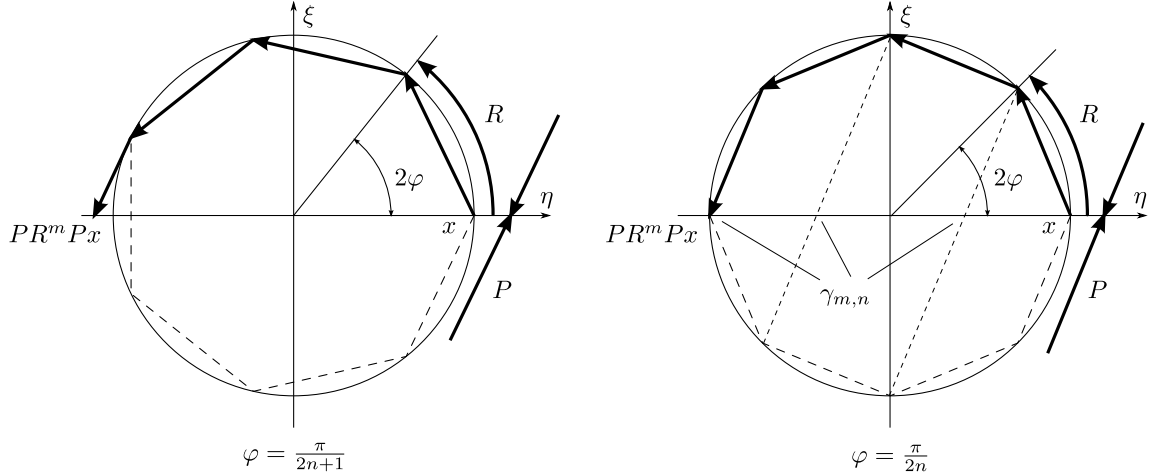}
\caption{Iterations of a point $x$
under the action of the map $PR^{m}P$}\label{semialg2}
\end{figure}

Both lemmas follow from the equality $(PR^{m}P)(\varphi)=\frac{\cos
[(2m+1)\varphi]}{\cos\varphi}P(\varphi)$ $(m=0,1, \ldots)$ whose validity is
justified by direct calculations, see Fig.~\ref{semialg2}.

\begin{lemma}\label{L:A3} Let $\varphi=\pi/(2n)$, and $B_{i}=P(\varphi)$
or $B_{i}=R(\varphi)$ for $1\le i\le m$. Then
\begin{equation}\label{E:A3}
B_{m}B_{m-1}\cdots B_{1}=\alpha R^{q}(\varphi)P^{r}(\varphi)R^{s}(\varphi),
\end{equation}
where $|\alpha |\le 1$, integers $q$ and $s$ are non-negative, $r=0,1$.
\end{lemma}

Prove the lemma by induction. For $m=1$ the assertion of the lemma is
evident; suppose that it is valid for $m=k-1\ge 1$. Then for $m=k$ the matrix
$A=B_{m}B_{m-1}\cdots B_{1}$ can be represented as
$A=B_{m}\tilde{\alpha}R^{\tilde{q}}P^{\tilde{r}}R^{\tilde{s}}_{m}$, where
$|\tilde{\alpha}|\le 1$, $\tilde{r}=0$ or $\tilde{r}=1$, $R=R(\varphi)$,
$P=P(\varphi)$.

If $B_{m}=R(\varphi)$ then
$A=\tilde{\alpha}R^{\tilde{q}+1}P^{\tilde{r}}R^{\tilde{s}}$, and for the
matrix $A$ the representation \eqref{E:A3} holds in which $\alpha=
\tilde{\alpha}$, $q=\tilde{q}+1$, $r=\tilde{r}$, $s=\tilde{s}$.

If $B_{m}=P(\varphi)$ and $\tilde{r}=0$ then
$A=\tilde{\alpha}PR^{\tilde{q}+\tilde{s}}$, and for the matrix $A$ the
representation \eqref{E:A3} holds in which $\alpha=\tilde{\alpha}$, $q=0$,
$r=1$, $s=\tilde{q}+\tilde{s}$.

If $B_{m}=P(\varphi)$ and $\tilde{r}=1$ then
$A=\tilde{\alpha}PR^{\tilde{q}}PR^{\tilde{s}}$. Here the factor
$PR^{\tilde{q}}P$ according to Lemma~\ref{L:A2} can be replaced by
$\lambda_{\tilde{q},n}P$. Then $A=\tilde{\alpha}\lambda_{\tilde{q},n}
PR^{\tilde{s}}$. Therefore, for the matrix $A$ the representation
\eqref{E:A3} holds in which $\alpha=\tilde{\alpha}\lambda_{\tilde{q},n}$,
$q=0$, $r=1$, $s=\tilde{s}$. In addition,
$|\alpha|\le|\tilde{\alpha}|\cdot|\lambda_{\tilde{q},n}|\le 1$ since
$|\tilde{\alpha}|\le 1$, $|\lambda_{\tilde{q},n}|\le 1$.

The inductive step is completed. Lemma~\ref{L:A3} is proved.

\begin{corollary*}
$|B_{m}B_{m-1}\cdots B_{1}|\le \left|P\left(\frac{\pi}{2n}\right)\right|$.
\end{corollary*}

The proof of the corollary immediates from the representation \eqref{E:A3}
and unitarity of the rotation matrix $R(\varphi)$.

\begin{lemma}\label{L:A4}
Let $t_{n}=\sin\frac{\pi}{2n}$. Then $\mathfrak{A}(t_{n})\in E(2,2)$.
\end{lemma}

Proof. Let $\{A(k)\}$ be a sequence of matrices from $\mathfrak{A}(t_{n})$.
Then, for each $k$, one of two equalities $A(k)=G(t_{n})$ or $A(k)= H(t_{n})$
holds. By \eqref{E:A1} $G(t_{n})=\mu_{n}P\left(\frac{\pi}{2n}\right)$,
$H(t_{n})=\mu_{n}R\left(\frac{\pi}{2n}\right)$, where
$\mu_{n}=1-\left(\sin\frac{\pi}{2n}\right)^{4}$. Therefore, the product of
matrices $A(1)$, $A(2)$, \ldots, $A(k)$ can be represented in the form:
$A(k)A(k - 1)\cdots A(1)=\mu ^{k}_{n}B_{k}B_{k-1}\cdots B_{1}$, where
$B_{i}=P\left(\frac{\pi}{2n}\right)$ or $B_{i}=R\left(\frac{\pi}{2n}\right)$.
Then, by Corollary from Lemma~\ref{L:A3}, $|A(k)A(k - 1)\cdots A(1)|\le \mu
^{k}_{n}\left|P\left(\frac{\pi}{2n}\right)\right|$ which implies absolute
stability of the class of matrices $\mathfrak{A}(t_{n})$. Lemma \ref{L:A4} is
proved.

\begin{lemma}\label{L:A5}
Let $s_{n}=\sin\frac{\pi}{2n+1}$. Then $\mathfrak{A}(s_{n})\not\in S(2,2)$
for all sufficiently large $n$.
\end{lemma}

Proof. Clearly the lemma will be proved if, for each sufficiently large $n$,
there can be found a sequence of matrices $A(k)\in\mathfrak{F}(s_{n})$ such
that
\begin{equation}\label{E:A4}
|A(k_{i})A(k_{i} - 1)\cdots A(1)| \rightarrow\infty
\end{equation}
for some $k_{i}\to\infty$.

Define the sequence of matrices $A(k)$ as follows: $A[(n+2)i]=G(s_{n})$,
$A[(n+2)i+1]=\ldots=A[(n+2)i+n]=H(s_{n})$, $A[(n+2)i+n+1]=G(s_{n})$. Let us
set $k_{i}=(n+2)i+n+1$. Then
\[
A(k_{i})A(k_{i}-1)\cdots A(1)=
[G(s_{n})H^{n}(s_{n})G(s_{n})]^{i}.
\]
Since $G(s_{n})=\nu_{n}P\left(\frac{\pi}{2n+1}\right)$ and
$H(s_{n})=\nu_{n}Q\left(\frac{\pi}{2n+1}\right)$, where
$\nu_{n}=1-\left(\sin\frac{\pi}{2n+1}\right)^{4}$, then
\[
A(k_{i})A(k_{i}-1)\cdots
A(1)=\left[\nu^{n+2}_{n}P\left(\frac{\pi}{2n+1}\right)Q^{n}\left(\frac{\pi}{2n+1}\right)P\left(\frac{\pi}{2n+1}\right)\right]^{i}.
\]
Consequently, by Lemma~\ref{L:A1} $A(k_{i})A(k_{i}-1)\cdots
A(1)=\left(-\nu^{n+2}_{n}/\cos\frac{\pi}{2n+1}\right)^{i}P\left(\frac{\pi}{2n+1}\right)$.
Recall that $P\left(\frac{\pi}{2n+1}\right)$ is a projector and so
$\left|P\left(\frac{\pi}{2n+1}\right)\right| \ge 1$. Therefore,
\begin{equation}\label{E:A5}
|A(k_{i})A(k_{i}-1)\cdots A(1)|\ge \left|-
\frac{\nu^{n+2}_{n}}{\cos\frac{\pi}{2n+1}}\right|^{i}.
\end{equation}

Direct calculations show that
$\nu^{n+2}_{n}/\cos\frac{\pi}{2n+1}=1+\frac{\pi^{2}}{2(2n+1)^{2}}+o(n^{-2})$.
Henceforth, for sufficiently large values of $n$ the inequality
$\nu^{n+2}_{n}/\cos\frac{\pi}{2n+1}>1$ holds. From this and from \eqref{E:A5}
the relation \eqref{E:A4} follows. Lemma~\ref{L:A5} is proved.

Let us complete the proof of the theorem. Since by Lemmas~\ref{L:A4} and
\ref{L:A5} $\mathfrak{A}(t_{n})\in E(2,2)$, $\mathfrak{A}(s_{n})\not\in
S(2,2)$ and $E(M,N)\subseteq U\subseteq S(M,N)$ then $\mathfrak{A}(t_{n})\in
U$, $\mathfrak{A}(s_{n})\not\in U$. But because of the points $t_{n}$ and
$s_{n}$ interleave between each other then the set $U\cap W$ contains
infinitely many different components of connectedness (different classes
$\mathfrak{A}(t_{n})$ belong to different components of connectedness).
Therefore, the set $U\cap W$ by the theorem of Whitney \cite{Milnor:e} is not
semialgebraic. But since $W$ is an algebraic set then the set $U$ is not
semialgebraic. Theorem~\ref{T:1} is proved.

2. \emph{Proof of Theorem~\ref{T:3}.} In one side the assertion of
Theorem~\ref{T:3} is obvious: absolute stability and absolute exponential
stability of the system \eqref{E:main} with respect to the class
$\mathfrak{A}=\{A_{1}, A_{2},\ldots,A_{M}\}$ immediately follow from
inequalities \eqref{E:ineq1} and \eqref{E:ineqq}.

Let us show that absolute exponential stability of the system \eqref{E:main}
with respect to the class of matrices $\mathfrak{A}=\{A_{1},
A_{2},\ldots,A_{M}\}$ implies \eqref{E:ineqq}. Let, for some $q<1$, the
relation \eqref{E:estim-exp} be valid. Set $\varkappa_{0}(x)=|x|$,
$\varkappa_{n}(x)=q^{-n}\max|B_{1}B_{2}\cdots B_{n}x|$ $(n>1)$, where
$|\cdot|$ is the Euclidean norm on $R^{N}$, and the maximum is taken over all
possible collections of the matrices $B_{1},B_{2},\ldots,
B_{n}\in\mathfrak{A}$. Define the norm $\|\cdot\|$ as follows: $\|x\| =
\sup_{n\ge0}\varkappa_{n}(x)$.

The function $\|x\|$ is semiadditive and due to \eqref{E:estim-exp} it
satisfies the relations $|x|\le\|x\|\le\max\{1,c\}\,|x|$. Henceforth,
$\|x\|=0$ if and only if $x=0$. Other properties of a norm are obvious for
$\|\cdot\|$.

Let us justify inequalities \eqref{E:ineqq}. Clearly, for any
$k=1,2,\ldots,M$ and $n=0,1,\ldots$ the estimate $\varkappa_{n}(A_{k}x)\le q
\varkappa_{n+1}(x)$ is valid. Therefore,	
\[
\|A_{k}x\|=\sup_{n\ge0}\varkappa_{n}(A_{k}x)\le q
\sup_{n\ge0}\varkappa_{n+1}(x)=\sup_{n\ge1}\varkappa_{n}(x)\le q\|x\|.
\]
From here $\|A_{k}\|\le q$ for $k=1,2,\ldots,M$. Inequalities \eqref{E:ineqq}
are proved.

Construction of the norm $\|\cdot\|$ and the proof of inequalities
\eqref{E:ineq1} in the case of absolute stability of the system
\eqref{E:main} are carried out similarly. Theorem~\ref{T:3} is proved.

  \providecommand{\bbljan}[0]{January} \providecommand{\bblfeb}[0]{February}
  \providecommand{\bblmar}[0]{March} \providecommand{\bblapr}[0]{April}
  \providecommand{\bblmay}[0]{May} \providecommand{\bbljun}[0]{June}
  \providecommand{\bbljul}[0]{July} \providecommand{\bblaug}[0]{August}
  \providecommand{\bblsep}[0]{September} \providecommand{\bbloct}[0]{October}
  \providecommand{\bblnov}[0]{November} \providecommand{\bbldec}[0]{December}

\section*{Post Scriptum}

In the foregoing text, a few misprints in the proof of Theorem~\ref{T:1}
occurred in the original journal version of the article \cite{Koz:AiT90:6:e}
were corrected, and two figures were added. The improved text was included in
the monograph \cite{AKKK:92:e}. Generalization of the presented results can
be found in \cite{Koz:AiT03:9:e}.

\end{document}